\documentclass[a4paper,12pt]{article}

\usepackage{amsmath, amssymb, amsthm}
\usepackage[dvips,final]{epsfig}
\usepackage[dvips,final]{graphicx}

\usepackage{color}

\definecolor{vio}{rgb}{0.5,0,0.5}
\definecolor{gre}{rgb}{0.1,0.7,0}
\definecolor{ora}{rgb}{0.8,0.2,0.1}
\definecolor{cya}{rgb}{0,.7,.1}

\def\red#1{\textcolor{black}{#1}}
\def\blu#1{\textcolor{black}{#1}}
\def\vio#1{\textcolor{black}{#1}}

\def\cyan#1{\textcolor{black}{#1}}

\usepackage{titlesec}
\titleformat{\section}{\bfseries}{\thesection}{1em}{}
\titleformat{\subsection}{\itshape}{\thesubsection}{1em}{}

\numberwithin{equation}{section}

\def\expe{\,\mathrm{e}}
\def\ve{\varepsilon}
\def\vp{\varphi}
\def\vr{\varrho}
\def\dd{\,\mathrm{d}}
\def\dive{\mathrm{\,div\,}}
\def\dist{\mathrm{dist}}

\def\real{\mathbb{R}}
\def\nat{\mathbb{N}}

\def \vect#1#2{\left(\!
\begin{array}{c} #1\\#2
\end{array}\! \right)}

\def \scal#1{\left\langle #1 \right\rangle}

\def\FF{\mathcal{F}}

\def\io{\int_{\Omega}}
\def\ipo{\int_{\partial\Omega}}

\def\om{^{(m)}}

\def\for{\mbox{ for }}
\def\Dom{\mathrm{Dom\,}}
\def\Span{\mathrm{Span}}
\def\Int{\mathrm{Int\,}}
\def\sign{\mathrm{sign\,}}
\def\supess{\mathop{\mbox{sup\,ess}}}

\def\be{\begin{equation}\label}
\def\ee{\end{equation}}

\def\ber{\begin{eqnarray}}
\def\eer{\end{eqnarray}}

\def\bers{\begin{eqnarray*}}
\def\eers{\end{eqnarray*}}

\newfont{\ctv}{msam10}

\newcommand{\bbox}{\mbox{\ctv \symbol{4}}}
\def\QED{{${}\hfill\bbox$}}
\newenvironment{pf}[1]{\par\vskip1mm{\noindent\it #1.}\ }{\QED\par
\vskip2mm}

\def\bpf{\begin{pf}}
\def\epf{\end{pf}}

\newtheorem{theorem}{Theorem}[section]

\newtheorem{coro}[theorem]{Corollary}
\newtheorem{hypo}[theorem]{Hypothesis}
\newtheorem{propo}[theorem]{Proposition}

\begin{document}

\title{Analysis of a tumor model as a multicomponent deformable porous medium}

\author{Pavel Krej\v c\'i
\thanks{Faculty of Civil Engineering, Czech Technical University, Th\'akurova 7, CZ-16629
Praha 6, Czech Republic, E-mail: {\tt Pavel.Krejci@cvut.cz}}
\and Elisabetta Rocca\thanks{Dipartimento di Matematica ``F. Casorati'', Universit\`a degli Studi di Pavia, and I.M.A.T.I. - CNR,
via Ferrata 5, Pavia I-27100, Italy,
E-mail: {\tt elisabetta.rocca@unipv.it}}
\and J\"urgen Sprekels\thanks{Weierstrass Institute for Applied Analysis and Stochastics, Mohrenstr.~39, D-10117 Berlin, Germany, E-mail: {\tt sprekels@wias-berlin.de}}}

\maketitle

\begin{abstract}
We propose a diffuse interface model to describe tumor as a multicomponent deformable porous medium. We include mechanical effects in the model by coupling the mass balance equations for the tumor species and the nutrient dynamics to a mechanical equilibrium equation with phase-dependent elasticity coefficients. The resulting PDE system couples two Cahn-Hilliard type equations for the tumor phase and  the healthy phase with a PDE linking the evolution of the interstitial fluid to the pressure of the system, a reaction-diffusion type equation for the nutrient proportion, and a quasistatic momentum balance. We prove here that the corresponding initial-boundary value problem has a solution in appropriate function spaces.
\end{abstract}

{\bf Key words:} tumor model, porous medium, diffuse interface model, Cahn-Hilliard equation, reaction-diffusion equation.


\section*{Introduction}\label{int}

Tumor growth is  nowadays one of the most active area of scientific research especially due to the impact on the quality of life for cancer patients. Starting with the seminal work of Burton \cite{Burton} and Greenspan \cite{Greenspan}, many mathematical models have been proposed to describe the complex biological and chemical processes that occur in tumor growth
with the aim of better understanding and ultimately controlling the behavior of cancer cells.
In recent years, there has been a growing interest in the mathematical modelling of cancer, see for example \cite{AgostiEtAl,Araujo,Bellomo,Byrne,Cristini,Fasano,Friedman}.
Mathematical models for tumor growth may have different analytical features: in the present work
we are focusing on the subclass of  continuum models, namely diffuse interface models. There are various ways to model the interaction between the tumor and surrounding host tissue. A classical approach is to represent the interfaces between the tumor and healthy tissues as idealized surfaces of zero thickness, leading to a sharp interface description that differentiates the tumor and the surrounding host tissue cell-by-cell. These sharp interface models are often difficult to analyze mathematically, and may fail when the interface undergoes a topological change.  Metastasis, which is the spreading of cancer to other parts of the body, is one important example of a change of topology. In  such an event, the interface can no longer be represented as a mathematical surface, and thus the sharp interface models do not properly describe the reality any more. 

On the other hand, diffuse interface models consider the interface between the tumor and the healthy tissues as a layer of non-infinitesimal thickness in which tumor and healthy cells can coexist. The main advantage of this approach is that the mathematical description is less sensitive to topological changes.  This is the reason why recent efforts in the mathematical modeling of tumor growth have been mostly focused on diffuse interface models, see for example \cite{CLLW,Cristini,Frieboes,GLNS,GLSS,Hawkins,Oden,Wise}, and their numerical simulations demonstrating complex changes in tumor morphologies
due to mechanical stresses and interactions with chemical species such as nutrients or toxic agents. 
Regarding the recent literature on the mathematical analysis of diffuse interface models for tumor growth we can further refer to \cite{CGH, CGRS1, CGRS2, Dai, FGR, FLR,  GLDirichlet, GLNeumann}
as mathematical references for Cahn-Hilliard-type models and 
\cite{BCG,GLDarcy,JWZ,LTZ} 
for models also including a transport effect described by Darcy's law.

A further class of diffuse interface models that also include chemotaxis and transport effects has been subsequently introduced (cf.~\cite{GLNS, GLSS}); moreover in some cases the sharp interface limits of such models have been investigated generally by using formal asymptotic methods (cf.~\cite{MR, RS}).

Including mechanics in the model is clearly an important issue that has been discussed in several modeling papers, but that has been very poorly studied analytically. Hence, the main aim of this paper is to find a compromise between the applications and the rigorous analysis of the resulting PDE system: we would like to introduce here an application-significant model which is tractable also analytically. 
Regarding the existing literature on this subject, we can quote the paper \cite{Sciume}, where, using multiphase porous media mechanics,
the authors represent a growing tumor as a multiphase medium containing an
extracellular matrix, tumor and host cells, and interstitial liquid. Numerical simulations are also
performed that characterize the process of cancer growth in terms of the initial tumor-to-healthy cell
density ratio, nutrient concentration, mechanical strain, cell adhesion, and geometry. However, referring to \cite{Tosin} for more details on this topic, we mention here that many models in the literature are based on the assumption that the tumor mass presents a particular geometry, the so-called spheroid, and in that case the models mainly focus on the evolution of the external radius of the spheroid. The resulting mathematical problem is an integro-differential free boundary problem, which has been proved to have solutions (cf.~\cite{Bueno,FrieRei}) and to predict the evolution of the system. Variants of this approach have been then considered, e.g. in \cite{CuiFrie} differentiating between viable cells and the necrotic core. Further extensions of the model introduced in \cite{Tosin} can be found in \cite{TosinPre}. 

Very recently in \cite{GLS} a new model for tumor growth dynamics including mechanical effects has been introduced in order to generalize the previous works \cite{Lima1, Lima2}  with the goal to take into account cell-cell adhesion effects  with the help of a Ginzburg-Landau type energy. In their model an equation of Cahn-Hilliard type is then coupled to the system of linear elasticity and a reaction-diffusion equation for a nutrient concentration and several questions regarding well-posedness and regularity of solutions are investigated.

In this paper, following the approach of \cite{Tosin}, we introduce a diffuse interface multicomponent model for tumor growth, where we include mechanics in the model, assuming that the tumor is a porous medium. In \cite{Tosin} the tumor is regarded as a mixture of various interacting components (cells and extracellular material)  whose evolution is ruled by  coupled mass and momentum balances. The cells usually are subdivided into subpopulations of proliferating, quiescent and necrotic cells (cf., e.\,g., \cite{CLLW, Cristini}) and the interactions between species are determined by the availability of some nutrients. 
Here, we restrict to the case where we distinguish only  healthy and tumor cells, even if we could, without affecting the analysis, treat the case where we differentiate also between necrotic and proliferating tumor cells. Hence, we represent the tumor as a porous medium consisting of three phases: healthy tissue $\vp_1$, tumor tissue $\vp_2$, and interstitial fluid $\vp_0$ satisfying proper mass balance equations including mass source terms depending on the nutrient variable $\vr$. The nutrient satisfies a reaction-diffusion equation nonlinearly coupled with the tumor and healthy tissue phases by a coefficient characterizing the different consumption rates of the nutrient by the different cell types. 
We couple the phases and nutrient dynamics with a mechanical equilibrium equation. This relation is further coupled with the phase dynamics through the elasticity modulus depending on the proportion between healthy and tumor phases. We refer to \cite{dk} for a mathematical model of a multicomponent flow in deformable porous media from which we take inspiration. The mass balance relations are derived from a free energy functional \cyan{${\mathcal F} = {\mathcal F}(\vp_0,\vp_1, \vp_2, w, \vr)$} which, in the domain $\Omega$ where the evolution takes place, can be written as 
\red{
\be{pote}
{\mathcal F} =\int_\Omega\left( \hat F(\vp_0 - w)+\frac{|\nabla\vp_1|^2}{2}+\frac{|\nabla\vp_2|^2}{2}+(\psi+g)(\vp_1,\vp_2)+\frac{\vr^2}{2} + \cyan{\frac{E(\vp_1,\vp_2)w^2}{2}}\right) \dd x 
\ee
where $w$ denotes the volume difference with respect to the referential state, $E$ is the elasticity modulus of the tissue, and $\hat F$ is a suitable nonnegative function defined below in \eqref{hatf}.} The sum $\psi+g$ represents the interaction potential \blu{of a typically double-well character} between tumor and healthy phases, with dominant component $\psi$ which is convex with bounded domain, while $g$ is its smooth non-convex perturbation. The quantity $\vr$ represents the mass content of the nutrient. Notice that the gradient terms in the free energy are due to the modeling assumption that the interface between healthy and tumor phases is diffuse (we take the parameters in front of the gradients equal to 1 here for simplicity, but, in practice, they determine the thickness of the interface and have to be chosen properly). The quantities $\vp_0, \vp_1, \vp_2$ are relative mass contents, so that only their nonnegative values are meaningful. We also assume that all the other substances present in the system are of negligible mass, that is, the identity $\vp_0 + \vp_1 + \vp_2 = 1$ is to be satisfied as part of the problem. Hence, we choose the domain of $\psi$ to be included in the set $\Theta:=\{(\vp_1,\vp_2)\in \real^2: \vp_1\ge 0,\ \vp_2\ge 0,\ \vp_1+\vp_2 \le 1\}$. Classically, $\psi$ can be taken as the indicator function of $\Theta$ or a logarithmic type potential (cf.~\cite{FLRS}).

Under proper assumptions on the data, we prove a result of existence of weak solutions for the resulting PDE system, that we will introduce in the next Session~\ref{mod}, 
coupled with suitable initial and boundary conditions. The PDEs consist of two Cahn-Hilliard type equations for the tumor phase and  the healthy phase with a PDE linking the evolution of the interstitial fluid to the pressure of the system, a reaction-diffusion type equation for the nutrient proportion and the momentum balance.  The technique of the proof is based on a regularization of the system, where, in particular, the non-smooth potential $\psi$ is replaced by his Yosida approximation $\psi_\ve$. Then, we prove existence of the approximated problem by means od a Faedo-Galerkin scheme and we pass to the limit by proving suitable uniform (in $\ve$) a-priori estimates and applying monotonicity and compactness arguments. A key point in the estimates consists in proving that the mean value of the phases are in the interior of the domain $\Theta$ of $\psi$, which in turns leads to the estimate of the mean value of the corresponding chemical potentials in the two Cahn-Hilliard type equations (cf.~\cite{ckrs, FLRS}). The uniqueness could be proved only in very particular situations, for example for smooth potentials $\psi$ satisfying suitable growth conditions and under some restrictions on the interaction coefficients in the Cahn-Hilliard type equations for the phase. We prefer to leave this argument for further studies on the model. 

{\bf Plan of the paper.} In the next Session~\ref{mod} we introduce the model deduced from the modeling hypothesis of \cite{Tosin} . In Session~\ref{sta} we state the mathematical problem and the main results of the paper concerning the existence of suitable weak solutions for the corresponding PDE system. The proof  relies on the passage to the limit (in Section~\ref{lim}) in a regularized problem, whose well posedness is obtained in Section~\ref{gal}.


\section{Modeling} \label{mod}

We follow the modeling hypotheses of \cite{Tosin} and represent the tumor as a porous medium consisting of three phases: healthy tissue, tumor tissue, and interstitial fluid. We choose the Lagrangian formalism, and assume that the evolution of the system takes place in a bounded domain $\Omega \subset \real^3$ with Lipschitzian boundary.

The state of the system is described by the following scalar quantities:
\begin{itemize}
\item[$\vp_0$:] Relative mass content of the interstitial fluid
\item[$\vp_1$:] Relative mass content of the healthy tissue
\item[$\vp_2$:] Relative mass content of the tumor tissue
\item[$\mu_1$:] Chemical potential controlling the growth of the healthy tissue
\item[$\mu_2$:] Chemical potential controlling the growth of the tumor tissue
\item[$p$:] Fluid pressure
\item[$w$:] Volume difference with respect to the referential state
\item[$\vr$:] Mass content of the nutrients
\end{itemize}
We consider the following evolution system in a given time interval $(0,T)$:
\begin{align}
\label{e1}
&\dot\vp_i+\sum_{j=0}^2 c_{ij} {\dive{\xi_j}}=S_i,\quad i=0,1,2,\\
\label{e2}
&\dot\vr+\dive{\zeta} + A(\vp_1,\vp_2)\, \vr = 0, \quad {\zeta}=-D\nabla \vr,\\
\label{e3}
&\nu\dot{w}+ E(\vp_1,\vp_2)\, w-p=\frac{1}{|\Omega|}\int_\Omega\big( E(\vp_1,\vp_2)\,w-p\big)\dd x , \\
\label{e4}
&\vect{\mu_1}{\mu_2}\in -\vect{\Delta\vp_1}{\Delta\vp_2}+\partial\psi(\vp_1, \vp_2) + \nabla_{\!\vp} g(\vp_1, \vp_2) \cyan{+ \nabla_{\!\vp} E(\vp_1, \vp_2) \frac{w^2}{2}},\\
\label{e5}
&S_0=-\gamma(\vr)\,\bar \vp_0(1-\vp_0), \quad S_1=\gamma (\vr)\,\bar\vp_0\vp_1, \quad S_2=\gamma(\vr)\,\bar\vp_0\vp_2,\\
\label{e6}
&\xi_j=-\nabla\mu_j\quad j=0,1,2, \\
\label{e7}
&\mu_0=p,\quad  w=\vp_0-f(p), 
\end{align}
where the dot denotes the derivative with respect to $t \in (0,T)$, $\partial\psi$ is the subdifferential of a convex potential $\psi$, $g$ is a smooth bounded possibly non-convex \blu{\fbox{removed}} perturbation of $\psi$, $\nabla$ is the gradient with respect to the space variable $x = (x_1, x_2, x_3)$, $\nabla_{\!\vp}$ is the gradient with respect to $\vp = (\vp_1, \vp_2)$, $\Delta$ is the Laplace operator, \red{$f$ is an empirical increasing function of the pressure,} and $\xi_j, \zeta$ are fluxes of the components $\vp_j, \vr$, respectively.

\red{The physical meaning of \eqref{e7} is the following: At constant volume $w$, the pressure $p$ increases when the  fluid content $\vp_0$ increases. Similarly, at constant pressure, the volume increases when the fluid content increases, and at constant fluid content, the pressure increases when the volume decreases.} 

The above system is coupled with initial and boundary conditions 
\begin{align}\label{ini1}
&\vp_i(0)=\vp_{}^0 \for i=1,2, \quad  w(0)=w^0, \quad \vr(0)=\vr^{0} \quad  \mbox{in }\ \Omega, \\ \label{bc1}
&\nabla\vp_i\cdot n=0\ \for i=1,2,\quad \xi_i\cdot n=0\ \for i=0,1,2,\quad {\zeta}\cdot n=\kappa (\vr-\vr^*)\ \ \mbox{on }\partial\Omega \times (0,T),
\end{align}
where $n = n(x)$ is the unit outward normal vector at the point $x \in \partial\Omega$.

In \eqref{e5} as well as in what follows, for a generic function $v \in L^1(\Omega\times (0,T))$ we denote by
\be{mean}
\bar v(t) = \frac{1}{|\Omega|}\io v(x,t) \dd x
\ee
for $t \in (0,T)$ the mean value of $v$ over $\Omega$.

\cyan{Equations} \red{\eqref{e4} and \eqref{e7} can be derived from the potential \eqref{pote} according to a standard ``Cahn-Hilliard'' theory, namely,
\be{cahi}
\mu_i \in \partial_{\vp_i} \FF(\vp_0,\vp_1, \vp_2, w, \vr), \quad i=0,1,2,
\ee
where $\partial_{\vp_i}$ denotes a suitable (e.\,g., Clarke) concept of subdifferential, provided we set
\be{hatf}
\hat F (z) = \int_{z_0}^z f^{-1}(s) \dd s
\ee
with $z_0$ from Hypothesis \ref{h1}\,(iv). Eqs.~\eqref{e1} then} represent the mass balance for the three components $\vp_0,\vp_1, \vp_2$ of the system, where $c_{ij}$ are the constant interaction coefficients. \red{The source terms $S_i$ are given by \eqref{e5}, where $\gamma(\rho)$ represents the speed of the growth rate depending on the nutrient concentration $\vr$.  Eq.~\eqref{e2} is a diffusion equation describing the mass balance for the nutrient concentration} with a constant positive diffusion coefficient $D>0$ and with a nonnegative coefficient $A$ depending on $\vp_1, \vp_2$ and characterizing the different consumption rates of the nutrient by the different cell types. The coefficient $\kappa>0$ in the boundary condition \eqref{bc1} for $\zeta$ is the diffusivity of the boundary for the nutrients, and $\vr^*$ is the (given) nutrient concentration outside the domain. Eq.~\eqref{e3} is the mechanical equilibrium equation with constant viscosity coefficient $\nu>0$ and with positive elasticity modulus $E(\vp_1,\vp_2)$ of the tissue which can be different for different proportions of $\vp_1$ and $\vp_2$. The constitutive functions $A, E, f, \gamma$, the convex potential $\psi$, the interaction constants, and the initial and boundary conditions satisfy Hypothesis \ref{h1} below.

\cyan{The system \eqref{e1}--\eqref{e7} is thermodynamically consistent.} \red{Indeed, the balance between the power supplied to the system $\io\sum_{i=0}^2 \dot\vp_i \mu_i \dd x$ and the potential increment $\dot \FF$ formally gives
\be{ther}
\begin{aligned}
\io\sum_{i=0}^2 \dot\vp_i \mu_i \dd x - \cyan{\frac{\dd}{\dd t}}\FF(\vp_0,\vp_1, \vp_2, w, \vr) &= \io (\dot w p - \cyan{E(\vp_1,\vp_2) \dot w w} - \dot\vr \vr)\dd x\\
&= \io (\nu |\dot w|^2 + D|\nabla\vr|^2 + A(\vp_1, \vp_2) \vr^2)\dd x \ge 0, 
\end{aligned}
\ee
which shows that the dissipation rate is positive during the process.}


\section{Statement of the problem} \label{sta}

The quantities $\vp_0, \vp_1, \vp_2$ are relative mass contents, so that only their nonnegative values are meaningful. We also assume that all the other substances present in the system are of negligible mass, that is, the identity $\vp_0 + \vp_1 + \vp_2 = 1$ is to be satisfied as part of the problem. The convex functional $\psi$ has to be chosen in such a way that the closure $\overline{\Dom \psi}$ of its domain $\Dom \psi$ is the set
\be{do1}
\overline{\Dom \psi} = \Theta := \{\vp = (\vp_1, \vp_2)\in \real^2: \vp_1\ge 0,\ \vp_2\ge 0,\ \vp_1+\vp_2 \le 1\},
\ee
and for $\delta\in (0,1 - (1/\sqrt{2}))$ we define
\be{do3}
\Theta_\delta := \{\vp \in \Int\Theta : \dist(\vp, \partial\Theta) \ge \delta\}.
\ee

Let us first specify the hypothesis about the data of the problem.

\begin{hypo}\label{h1}
We fix a constant $K\ge 1$ and assume the following hypothesis to hold.
\begin{itemize}
\item[{\rm (i)}] $\displaystyle{\sum_{i=0}^2 c_{ij}} = 0$ for all $j=0,1,2$, $\displaystyle{\sum_{j=0}^2 c_{ij}} = 0$ for all $i=0,1,2$, and there exists $\hat c > 0$ such that 
$\displaystyle{-\sum_{i\ne j} c_{ij} |\xi_i - \xi_j|^2 \ge \hat c\big(|\xi_1 - \xi_0|^2+|\xi_2 - \xi_0|^2\big)}$ for all $\xi_0, \xi_1, \xi_2 \in \real^3$;
\item[{\rm (ii)}] $E, A : \real^2 \to [0,K]$ are Lipschitz continuous functions;
\item[{\rm (iii)}] $\gamma : \real \to [-K,K]$ is a continuously differentiable function, $|\gamma'(\vr)| \le K$ for all $\rho\in \real$;
\item[{\rm (iv)}] $f : \real \to \real$ is a continuously differentiable function, \vio{$f_1 \ge f'(p) \ge f_0$ for some $f_1 > f_0 > 0$} and all $p\in \real$, \red{$z_0 = f(0)$};
\item[{\rm (v)}] $\psi: \real^2 \to [0,+\infty]$ is a proper convex lower semicontinuous function satisfying \eqref{do1}.
We further assume that there exist positive constants $\delta, b', c', r'$ such that putting $\delta_T = \delta\expe^{-KT-2}$ \blu{\fbox{removed}} the following implications hold:
\begin{itemize}
\item[{\rm (v1)}] $\dist(\hat\vp,\Theta_{\delta_T}) \le \delta_T/2 \ \Longrightarrow \ |\hat\xi| \le b' \quad \forall \hat\xi \in \partial \psi(\hat\vp)$;
\item[{\rm (v2)}] $\hat\vp \in \Theta_{\delta_T},\ |\vp-\hat\vp|\ge \delta_T/4, \ \vp \in \Theta \Longrightarrow\ r'|\xi - \hat\xi|\leq \scal{\xi-\hat\xi,\vp - \hat\vp}+ c'$
\item[] \qquad $\forall \xi \in \partial\psi(\vp), \quad \forall \hat\xi \in \partial\psi(\hat\vp)$;
\end{itemize}
\item[{\rm (vi)}] $g : \Theta \to \real$ is a given function of class $C^2$;
\item[{\rm (vii)}] $\vp_0^0, \vp_1^0, \vp_2^0, w^0, \vr^0 \in W^{1,2}(\Omega)\cap L^\infty(\Omega)$ are given initial \blu{data} such that $\overline{w^0}=0$, $(\overline{\vp_1^0}, \overline{\vp_2^0}) \in \Theta_\delta$ with $\delta$ from Hypothesis (v), $\vp_0^0(x)+\vp_1^0(x)+ \vp_2^0(x) =1$ for a.\,e. $x \in \Omega$;
\item[{\rm (viii)}] $\vr^* \in L^\infty((0,T)\times \partial\Omega)$ is a given function with $\dot\vr^* \in L^2((0,T)\times \partial\Omega)$.
\end{itemize}
\end{hypo}

Conditions (v1), (v2) need some comments. They slightly differ from those in \cite[Proposition 2.10]{ckrs}, but it is easy to check they are still satisfied if for example $\psi$ is the indicator function of the set $\Theta$. Indeed, (v1) holds trivially. To verify that (v2) holds, take any $\vp \in \Theta$ and $\xi \in \partial\psi(\vp)$. We first notice that $\hat\xi = 0$, and
$$
\scal{\xi, \vp - v} \ge 0 \quad \forall v \in \Theta.
$$
We are done if $\xi = 0$. Otherwise,
$$
v = \hat\vp + \delta_T\frac{\xi}{|\xi|}
$$
is an admissible choice, and we obtain
$
\scal{\xi, \vp - \hat\vp} \ge \delta_T|\xi|,
$
which is precisely (v2) with $r'=\delta_T$ and $c'=0$.

In the proof, we extend the function $g$ to the whole $\real^2$ in such a way that
\be{cg}
C_g := \sup\{|g(\vp)|, |\nabla_{\!\vp}g(\vp)|, |\scal{\nabla_{\!\vp}g(\vp), \vp}|: \vp \in \real^2\} < \infty.
\ee

The main result of the paper reads as follows.

\begin{theorem}\label{t1}
Let Hypothesis \ref{h1} hold. Then the system \blu{\eqref{e1}--\eqref{bc1}} admits a solution with the regularity $\vp_i \in L^\infty((0,T)\times\Omega)$, $\nabla \vp_i \in L^\infty(0,T;L^2(\Omega))$, $\dot \vp_i \in L^2(0,T;W^{-1,2}(\Omega))$, $\mu_i, \nabla \mu_i \in L^2((0,T)\times\Omega)$ for $i=0,1,2$, $(\vp_1(x,t), \vp_2(x,t)) \in \Theta$ a.\,e., $\vp_0+\vp_1+ \vp_2 =1$ a.\,e., $w \in L^\infty((0,T)\times\Omega)$, $\dot w, \nabla w, \nabla\dot w \in L^\infty(0,T;L^2(\Omega))$, $\dot\vr \in L^2((0,T)\times\Omega)$, $\vr, \nabla \vr \in L^\infty(0,T;L^2(\Omega))$. The equations \eqref{e3}, \eqref{e5}--\eqref{e7} \blu{and the initial conditions \eqref{ini1}} are satisfied almost everywhere in $\Omega\times(0,T)$ \blu{and in $\Omega$, respectively, and the relations} \eqref{e1}--\eqref{e2} and \eqref{e4} are to be interpreted respectively as
\begin{align}
\label{e1u}
&\io \left(\dot \vp_i\, v_i + \sum_{j=0}^2 c_{ij} \scal{\nabla \mu_j, \nabla v_i}\right)\dd x = \io S_i\,v_i \dd x,\quad i=0,1,2,\\ \label{e2u}
&\io \big(\dot\vr \hat v+D \scal{\nabla\vr,\nabla \hat v} + A(\vp_1,\vp_2)\, \vr\hat v\big)\dd x + \kappa \ipo (\vr - \vr^*)\,\hat v \dd s(x) = 0,\\ \nonumber
&\cyan{\io\scal{\vect{\mu_1}{\mu_2} - \nabla_{\!\vp}g(\vp_1, \vp_2) - \nabla_{\!\vp} E(\vp_1, \vp_2)\frac{w^2}{2},\vect{v_1}{v_2}{-}\vect{\vp_1}{\vp_2}} \dd x}\\ \label{e4u}
&\quad -\io\big(\scal{\nabla\vp_1,\nabla(v_1{-}\vp_1)} + \scal{\nabla\vp_2,\nabla(v_2{-}\vp_2)}\big)\dd x \le \io\big(\psi(v_1, v_2) -\psi(\vp_1, \vp_2) \big)\dd x
\end{align}
for a.\,e. $t \in (0,T)$ and for all test functions $v_0, v_1, v_2, \hat v\in W^{1,2}(\Omega)$.
\end {theorem}

The proof of Theorem \ref{t1} is divided into several steps. We introduce a small regularizing parameter $\ve > 0$ and approximate
the convex potential $\psi$ by its Yosida approximation $\psi^\ve$ defined by the formula
\begin{equation}\label{yosi}
\psi^\ve(\vp) = \min_{z\in \real^2}\left\{\frac{1}{2\ve}|
\vp-z|^2+\psi(z)\right\}.
\end{equation}
Let us recall the main properties of the Yosida approximation, see
\cite{ae,barbu,brezis} for proofs.

\begin{propo}\label{yos}
The mapping $\psi^\ve:\real^2 \to [0,\infty)$ is convex and continuously differentiable, and the so-called {\em resolvent}
$J^\ve$ of $\partial \psi$, defined as
\begin{equation}\label{jep}
J^\ve = (I+\ve\,\partial\psi)^{-1},
\end{equation}
where $I$ is the identity, is non-expansive in $\real^2$. The mapping $\nabla_{\!\vp}\psi^\ve$ is monotone and Lipschitz continuous, and has for every $\vp \in \real^2$ the properties
\begin{align}\label{yos1}
&\nabla_{\!\vp}\psi^\ve(\vp) = \frac{1}{\ve}(\vp-J^\ve \vp)
\in \partial\psi(J^\ve \vp) \quad\forall\ve>0,\\
\label{yos2}
&\vp\in \Dom\partial\psi \ \Longrightarrow \
\left\{
\begin{array}{l}
|\nabla_{\!\vp}\psi^\ve(\vp)- m(\partial\psi(\vp))|\to 0\\[1mm]
|\nabla_{\!\vp}\psi^\ve(\vp)|\nearrow |m(\partial\psi(\vp))|
\end{array}
\right.
\quad\hbox{as }\ \ve\searrow 0, \\
\label{yos3}
&\psi^\ve(\vp)=\frac{\ve}{2}|\nabla_{\!\vp}\psi^\ve(\vp)|^2+\psi(J^\ve
\vp) \quad\forall \ve>0,\\
\label{yos4}
&\psi^\ve(\vp)\nearrow\psi(\vp)\quad\hbox{as }\ \ve\searrow 0,
\end{align}
where $m(\partial\psi(\vp))$ is the element of $\partial\psi(\vp)$ with minimal norm.
\end{propo}

{}From \eqref{yos1}--\eqref{yos3} it follows that for every $\vp \in \real^2$ and every $\ve >0$ we have
\be{yos5}
\psi^\ve(\vp) = \frac{1}{2\ve}|\vp - J^\ve\vp|^2 + \psi(J^\ve\vp).
\ee
\blu{For every $\vp\in \real^2$ and $\ve>0$ we have $J^\ve(\vp) \in \Theta$, so that $\psi(J^\ve\vp) \ge 0 \ge |J^\ve\vp|^2 -1$. Furthermore, the Young inequality yields that
$$
2\scal{\vp,J^\ve\vp} \le \frac{1}{2\ve + 1} |\vp|^2 + (2\ve + 1) |J^\ve\vp|^2,
$$
and} we obtain
\be{lowps}
\psi^\ve(\vp) \ge \frac{1}{2\ve + 1}|\vp|^2 - 1 \quad \forall \vp \in \real^2.
\ee
We consider the following weak formulation of the regularized problem \blu{\eqref{e1}--\eqref{bc1}}
\begin{align}
\label{e1w}
&\io \left(\dot \vp_i\, v_i + \sum_{j=0}^2 c_{ij} \scal{\nabla \mu_j, \nabla v_i}\right)\dd x = \io S_i\,v_i \dd x,\quad i=0,1,2,\\ \label{e2w}
&\io \big(\dot\vr \hat v+D \scal{\nabla\vr,\nabla \hat v} + A(\vp_1,\vp_2)\, \vr\hat v\big)\dd x + \kappa \ipo (\vr - \vr^*)\,\hat v \dd s(x) = 0,\\
\label{e3w}
&\nu\dot{w} + E(\vp_1,\vp_2) w-\frac{p}{|\vp_0|{+}|\vp_1|{+}|\vp_2|} = \frac{1}{|\Omega|}\int_\Omega\left( E(\vp_1,\vp_2)\,w-\frac{p}{|\vp_0|{+}|\vp_1|{+}|\vp_2|}\right)\dd x, \\
\label{e4w}
&\vect{\mu_1}{\mu_2} = -\vect{\Delta\vp_1}{\Delta\vp_2}+\nabla_{\!\vp}\left(\psi^\ve(\vp_1, \vp_2) + g(\vp_1, \vp_2)\right) \cyan{+\nabla_{\!\vp} E(\vp_1,\vp_2) \frac{w^2}{2},}\\
\label{e5w}
&S_0=-Q\,(1-\vp_0), \ \ S_1=Q\,\vp_1, \ \ S_2=Q\, \vp_2,\\[2mm] \label{qr}
&Q =\frac{\gamma(\vr)\,\bar\vp_{0}}{(|\vp_0|+|\vp_1|+|\vp_2|)(|\bar\vp_0|+|\bar\vp_1|+|\bar\vp_2|)},\\[2mm]
\label{e6w}
&\xi_j=-\nabla\mu_j\quad j=0,1,2, \\
\label{e7w}
&\mu_0=p,\quad  w=\vp_0-f(p), 
\end{align}
for a.\,e. $t \in (0,T)$ and for all test functions $v_0, v_1, v_2, \hat v\in W^{1,2}(\Omega)$.

Assuming that \eqref{e1w}--\eqref{e7w}, \eqref{ini1} has a solution, choosing $v_0 = v_1 = v_2 = v$ in \eqref{e1w}, and summing up over $i=0,1,2$, we obtain formally from Hypothesis \ref{h1}\,(i) the identity
$$
\io \left(\sum_{i=0}^2\dot \vp_i\right) v \dd x= \io Q\left(\sum_{i=0}^2\vp_i-1\right)v\dd x
$$
for all $v \in W^{1,2}(\Omega)$. \red{Putting $y(x,t) = \sum_{i=0}^2\vp_i(x,t)-1$, we see that this is an identity of the form $\dot y(x,t) = Q(x,t) y(x,t)$ with initial condition $y(x,0) = 0$ according to Hypothesis \ref{h1}\,(vii). Hence,} still formally,
\be{ide}
\sum_{i=0}^2\vp_i(x,t) = 1
\ee
for all $x$ and $t$. In particular, the denominators in \eqref{e3w} and \eqref{qr} are greater or equal to one. We show below that in the limit $\ve \to 0$, all $\vp_i$ will be nonnegative, and the denominators will all be equal to $1$.


\section{Galerkin approximations}\label{gal}

We solve the problem \eqref{e1w}--\eqref{e7w}, \eqref{ini1} by Galerkin approximations.  We choose the orthonormal basis ${e_k: k\in \nat\cup \{0\}}$ in $L^2(\Omega)$ such that 
\[
-\Delta e_k=\lambda_k e_k \quad \hbox{ in }\Omega, \quad  \nabla e_k\cdot n =0\quad\hbox{ on }\partial\Omega\ \for k \in \nat\cup \{0\}, \quad \lambda_0=0,
\]
and for $m \in \nat$ we introduce the functions
\begin{align*}
&\vp_i^{(m)}(x,t)=\sum_{k=0}^m \tilde\vp_{ik}(t) e_k(x), \quad \mu_i^{(m)}=\sum_{k=0}^m \tilde\mu_{ik}(t) e_k(x) \ \for i=0,1,2,\\ 
&\vr^{(m)}(x,t)=\sum_{k=0}^m \tilde\vr_{k}(t) e_k(x),
\end{align*}
with time dependent coefficients $\tilde\vp_{ik}(t), \tilde\mu_{ik}(t), \tilde\vr_{k}(t)$ which are to be found as solutions of the ODE system for $k=0, 1, \dots, m$
\begin{align}
\label{e1m}
&\io \left(\dot \vp_i\om e_k  + \sum_{j=0}^2 c_{ij} \scal{\nabla \mu_j\om, \nabla e_k}\right)\dd x = \io S_i\om\,e_k \dd x,\quad i=0,1,2,\\ \label{e2m}
&\io \big(\dot\vr\om e_k+D \scal{\nabla\vr\om,\nabla e_k} + A(\vp_1\om,\vp_2\om)\, \vr\om e_k\big)\dd x + \kappa \ipo (\vr\om - \vr^*) e_k \dd s(x) = 0,\\ \nonumber
&\nu \dot{w}\om+ E(\vp_1\om,\vp_2\om)\, w\om- \frac{f^{-1}(\vp_0\om - w\om)}{|\vp_0\om|+|\vp_1\om|+|\vp_2\om|} \\ \label{e3m}
&\hspace{14mm} =\frac{1}{|\Omega|}\int_\Omega\left( E(\vp_1\om,\vp_2\om)\,w\om-\frac{f^{-1}(\vp_0\om - w\om)}{|\vp_0\om|+|\vp_1\om|+|\vp_2\om|}\right)\dd x ,
\end{align}
\begin{align} \label{e4m}
\mu_0\om &= P_m(f^{-1}(\vp_0\om - w\om)),\\
\label{e4mm}
\mu_i\om &= -\Delta\vp_i\om + P_m\left(\partial_i\psi^\ve(\vp_1\om, \vp_2\om)+\partial_i g(\vp_1\om, \vp_2\om) \cyan{+ \partial_i E(\vp_1\om, \vp_2\om) \frac{|w\om|^2}{2}}\right), \ i=1,2,\\
\label{e5m}
S_0\om &=-Q\om\,(1-\vp_0\om), \ \ S_1\om=Q\om\,\vp_1\om, \ \ S_2\om=Q\om\, \vp_2\om,\\[2mm] \label{qrm}
Q\om &= \frac{\gamma(\vr\om)\,\bar\vp_{0}\om}{(|\vp_0\om|+|\vp_1\om|+|\vp_2\om|)(|\bar\vp_0\om|+|\bar\vp_1\om|+|\bar\vp_2\om|)},
\end{align}
where $P_m: L^2(\Omega) \to H_m := \Span(e_0, \dots, e_m)$ is the orthogonal projection of $L^2(\Omega)$ onto $H_m$, \cyan{the symbol $\partial_i$ denotes the partial derivative $\partial/\partial \vp_i$ for $i=1,2$, and} $\bar\vp_{i}\om = (1/|\Omega|)\io\vp_{i}\om\dd x$. The initial conditions are
\be{inik}
\tilde\vp_{ik}(0)= \io \vp_{i}^0(x)\,e_k(x)\dd x, \ \ \tilde\vr_k(0)=\io\vr^{0}(x)\,e_k(x)\dd x,\ \ \vio{ w\om(x,0)= w^0(x).}
\ee
System \eqref{e1m}--\eqref{e2m} is a locally well-posed system of $4(m+1)$ differential equations of the first order for $4(m+1)$ scalar unknowns $\tilde \vr_k, \tilde\vp_{ik}$, $i=0,1,2$, $k=0,1,\dots m$, while it is convenient to interpret \eqref{e3m}--\eqref{e5m} as constitutive relations. We shall see below in Eq.~\eqref{es0} that the expressions in the denominators of \eqref{e3m} and \eqref{qrm} are greater or equal to $1$, hence the formulas are meaningful. In particular, since $f^{-1}$ is Lipschitz continuous by Hypothesis \ref{h1}\,(iv), the equation \eqref{e3m} defines a Lipschitz continuous solution operator $W:C([0,T];\real^{3(m+1)}) \to C^1([0,T]; W^{1,2}(\Omega))$ which with given functions $\tilde\vp_{ik}$, $i=0,1,2$, $k=0,1,\dots m$ associates the solution $w\om$ of \eqref{e3m}. The existence of a unique local solution of \eqref{e1m}--\eqref{e5m} is therefore guaranteed on a nondegenerate time interval $[0,T_m)$, $0< T_m \le T$.

In order to show that the solution \eqref{e1m}--\eqref{e5m} is global, we derive some estimates for the solution on the whole interval $[0,T_m)$.


\subsection{Estimates independent of $m$}\label{es01}

In the series of estimates which we derive in the formulas below, we denote by $C$ any positive constant independent of $m$ and $\ve$, and by $C^\ve$ any constant independent of $m$ and depending possibly on $\ve$. For simplicity, we denote by $|\cdot|_H$ the norm in $L^2(\Omega)$, and by $\|\cdot \|_V$ the norm in $W^{1,2}(\Omega)$.

We first handle Eq.~\eqref{e2m} which is easy. We multiply it by $\tilde \vr_k$ and sum up over $k=0, \dots, m$ to obtain
\be{esro}
\frac{\dd}{\dd t}\frac12 \io |\vr\om|^2\dd x + D\io |\nabla\vr\om|^2\dd x + \frac{\kappa}{2} \ipo |\vr\om|^2 \dd s(x) \le C.
\ee
We proceed similarly multiplying \eqref{e2m} by $\dot{\tilde\vr}_k$ and summing up over $k=0, \dots, m$ to obtain
that
\be{esro2}
\io |\dot\vr\om|^2\dd x + \frac{\dd}{\dd t}\left(D \io |\nabla\vr\om|^2\dd x + \kappa \ipo |\vr\om|^2 \dd s(x)\right) \le C\left(1+\io |\vr\om|^2\dd x\right),
\ee
hence,
\be{esro3}
\int_0^{T_m}\io |\dot\vr\om|^2\dd x\dd t + \supess_{t\in (0,T_m)}\left( \io |\nabla\vr\om|^2(t)\dd x + \ipo |\vr\om|^2(t) \dd s(x)\right) \le C.
\ee
We further \cyan{sum up Eqs.~\eqref{e1m} over $i=0,1,2$. From Hypothesis \ref{h1}\,(i) it follows that
\begin{align} \nonumber
\frac{\dd}{\dd t} \io(\vp_0\om {+} \vp_1\om {+}\vp_2\om)\,e_k \dd x &= \!\io(S_0\om {+} S_1\om {+}S_2\om)\,e_k \dd x\\ \label{es0a} 
& = \!\io Q\om(x,t)\big(\vp_0\om{+}\vp_1\om{+}\vp_2\om - 1\big)e_k \dd x
\end{align}
for all $k=0,\dots, m$. In terms of the functions
$$
y_k(t) = \io(\vp_0\om {+} \vp_1\om {+}\vp_2\om -1)\,e_k\dd x = \tilde \vp_{0,k} + \tilde \vp_{1,k} + \tilde \vp_{2,k} - \delta_{0,k},
$$
where $\delta_{0,k}$ is the Kronecker symbol, we rewrite \eqref{es0a} in the form
$$
\dot y_k(t) = \sum_{l=0}^m a_{kl}(t) y_l(t)
$$ 
with bounded coefficients $a_{kl}(t)$. This is a linear ODE system with zero initial conditions by Hypothesis \ref{h1}\,(vii), so that all functions $y_k(t)$ vanish in the whole interval of existence. Hence,}
\be{es0}
\vp_0\om(x,t) + \vp_1\om(x,t) +\vp_2\om(x,t) = 1
\ee
for all $(x,t) \in \Omega\times [0,T_m)$.

\cyan{Next, we prove that} $w\om$ as a solution of the ODE \eqref{e3m} admits, \cyan{together with its time derivative}, an $L^\infty$-bound independent of $m$ and $\ve$, namely
\be{es10}
\supess_{(x,t)\in \Omega\times (0,T_m)} (|w\om(x,t)| \cyan{+ |\dot w\om(x,t)|}) \le C.
\ee
Indeed, we first add to both the left-hand side and the right-hand side of \eqref{e3m} the term
$$
\frac{f^{-1}(\vp_0\om)}{|\vp_0\om|+|\vp_1\om|+|\vp_2\om|}
$$
which is bounded by Hypothesis \ref{h1}\,(iv) \cyan{and \eqref{es0}.} 
Note that by \eqref{e3m} for $k=0$ we have $\io \dot w\om \dd x =0$, hence $\io w\om \dd x =0$ by Hypothesis \ref{h1}\,(vii). Next, we multiply \eqref{e3m} by $w\om$, use the fact that the mean value of $w\om$ is zero and that \red{$(f^{-1}(\vp_0\om) - f^{-1}(\vp_0\om {-} w\om))\, w\om\ge 0$}, integrate over $\Omega$, and obtain
$$
\frac{\dd}{\dd t} \io |w\om(x,t)|^2 \dd x \le \red{C\io |w\om(x,t)| \dd x}
$$
for a.\,e. $t\in (0,T_m)$, hence, $\io |w\om(x,t)|^2 \dd x \le C$ for $t\in [0,T_m)$. In particular, the right-hand side of \eqref{e3m} is bounded independently of $m$ and $\ve$. We now repeat the same procedure multiplying \eqref{e3m} by $\sign w\om$ without integration over $\Omega$. \blu{Note that for a.\,e. $x \in \Omega$, the function $t \mapsto w\om(x,t)$ is absolutely continuous, so that $\dot w\om(x,t)\, \sign w\om(x,t) = \frac{\partial}{\partial t} |w\om(x,t)|$ a.\,e. Furthermore, $(f^{-1}(\vp_0\om) - f^{-1}(\vp_0\om {-} w\om))\,\sign w\om \ge 0$, and we get
$$
\nu \frac{\partial}{\partial t} |w\om(x,t)| + E(\vp_1\om,\vp_2\om)|w\om(x,t)| \le C \quad \mbox{a.\,e. in }  \Omega\times(0,T_m)
$$
with $|w\om (x,0)| \le C$. Integrating from $0$ to $t$} we \cyan{get a uniform upper bound for $|w\om(x,t)|$, and, by comparison in \eqref{e3m}, we} conclude that \eqref{es10} holds.

\cyan{Further estimates are} more delicate. 
We multiply the $(i,k)$-th equation of \eqref{e1m} by $\tilde \mu_{ik}$ and sum up over $i=0,1,2$ and $k=0,1,\dots, m$ to obtain
\be{es1}
\sum_{i=0}^2\io \dot \vp_i\om \mu_i\om \dd x + \sum_{i,j=0}^2 c_{ij} \io\scal{\nabla \mu_j\om, \nabla \mu_i\om}\dd x = \sum_{i=0}^2\io S_i\om\,\mu_i\om \dd x.
\ee
We treat the three integrals in \eqref{es1} separately. \cyan{We first define a reduced potential $\mathcal{F}^\ve_0 = \mathcal{F}^\ve_0(\vp_0, \vp_1, \vp_2, w)$ by the formula
\be{poteps}
\mathcal{F}^\ve_0 = \io\left(\hat F(\vp_0 - w) + \psi^\ve(\vp_1, \vp_2) + g(\vp_1, \vp_2) +E(\vp_1, \vp_2) \frac{w^2}{2} + \frac{|\nabla \vp_1|^2}{2} + \frac{|\nabla \vp_2|^2}{2}\right)\dd x
\ee
with $\hat F$ as in \eqref{hatf}. Then the first integral on the left-hand side of \eqref{es1} can be rewritten as
\begin{align} \nonumber
\sum_{i=0}^2\io \dot \vp_i\om \mu_i\om \dd x &= \frac{\dd}{\dd t} \mathcal{F}^\ve_0(\vp_0\om, \vp_1\om, \vp_2\om, w\om)\\ \nonumber
&\quad +\io \dot w\om (f^{-1}(\vp_0\om - w\om) - E(\vp_1\om,\vp_2\om)w\om)\dd x\\ \label{es2}
&\ge \frac{\dd}{\dd t} \mathcal{F}^\ve_0(\vp_0\om, \vp_1\om, \vp_2\om, w\om)  - C \io |\vp_0\om| \dd x
\end{align}
with a constant $C$ independent of $m$ and $\ve$ as a consequence of \eqref{es10}.}

To estimate the second integral in \eqref{es1}, we use the vector formula
$$
\scal{u,v} = -\frac12(|u-v|^2 - |u|^2 - |v|^2)
$$
to conclude, using Hypothesis \ref{h1}\,(i), that
\be{es4}
\sum_{i,j=0}^2 c_{ij}\scal{\nabla \mu_j\om, \nabla \mu_i\om}=
-\frac12\sum_{i\ne j} c_{ij}|\nabla \mu_j\om - \nabla \mu_i\om|^2 \ge \frac{\hat c}{2} \sum_{i=1}^2|\nabla \mu_i\om - \nabla \mu_0\om|^2.
\ee
Finally, the integral on the right-hand side of \eqref{es1} can be rewritten in the form
\begin{align}\nonumber
&\sum_{i=0}^2\io S_i\om\,\mu_i\om \dd x = - \sum_{i=1}^2\io S_i\om \Delta\vp_i\om \dd x + \io S_0\om\,f^{-1}(\vp_0\om - w\om) \dd x\\ \label{es5}
& \quad + \sum_{i=1}^2\io S_i\om\, P_m\left(\partial_i\psi^\ve(\vp_1\om, \vp_2\om) \cyan{+ \partial_i g(\vp_1\om, \vp_2\om)+ \partial_i E(\vp_1\om, \vp_2\om)\frac{|w\om|^2}{2}}\right) \dd x.
\end{align}
\cyan{The first term on the right-hand side of \eqref{es5} can be estimated using integration by parts as follows.
\be{es8}
-\sum_{i=1}^2\io S_i\om\,\Delta\vp_i\om \dd x = \sum_{i=1}^2\io \scal{\nabla S_i\om, \nabla\vp_i\om}\dd x \le C\left(|\nabla \vr\om|_H^2 {+} |\nabla \vp_1\om|_H^2 {+} |\nabla \vp_2\om|_H^2\right).
\ee
To estimate the remaining terms, notice that}
the function $Q\om$ defined in \eqref{qrm} is bounded in absolute value by the constant $K$ from Hypothesis \ref{h1}\,(iii), and also
\be{es6}
|S_i\om(x,t)| \le K \quad \mbox{for all } \ x\in \Omega, \ t\in [0,T_m), \ i=0,1,2\,.
\ee
By Proposition \ref{yos}, the gradient $\nabla_{\!\vp}\psi^\ve$ of $\psi^\ve$
is Lipschitz continuous with a constant depending on $\ve$. We thus obtain from \eqref{es5}\cyan{--\eqref{es6} that
\be{es7}
\sum_{i=0}^2\io S_i\om\,\mu_i\om \dd x \le C^\ve \left(1{+}|\vp_0\om|_H {+} |\vp_1\om|_H {+} |\vp_2\om|_H {+} |\nabla \vr\om|_H^2 {+} |\nabla \vp_1\om|_H^2 {+} |\nabla \vp_2\om|_H^2\right).
\ee
}
Combining \eqref{es1}--\eqref{es7} we thus obtain
\cyan{
\begin{align}\nonumber
&\frac{\dd}{\dd t} \cyan{\mathcal{F}_0^\ve(\vp_0\om, \vp_1\om, \vp_2\om, w)} + \io \left(\sum_{i=1}^2 |\nabla \mu_i\om - \nabla \mu_0\om|^2\right) \dd x \\ \label{es9}
&\le C^\ve\left(1 + |\nabla \vr\om|_H^2 + |\nabla \vp_1\om|_H^2 + |\nabla \vp_2\om|_H^2 +|\vp_0\om|_H^2 + |\vp_1\om|_H^2 + |\vp_2\om|_H^2\right).
\end{align}
By Hypothesis \ref{h1}\,(iv), we have $\hat F(z) \ge z^2/(4f_1) - C$ for all $z\in\real$ with some constant $C>0$. Hence, using \eqref{poteps}}, \eqref{lowps}, \eqref{esro3}, \eqref{es10}, and Gronwall's argument, we derive from \eqref{es9} the estimate
\begin{align}\nonumber
&\supess_{t\in (0,T_m)}\left(|\vp_0\om|_H^2(t) + |\vp_1\om|_H^2(t) + |\vp_2\om|_H^2(t) + |\nabla \vp_1\om|_H^2(t) + |\nabla \vp_2\om|_H^2(t)\right)\\ \label{es9a}
& \hspace{10mm} +
\int_0^{T_m} \io \left(\sum_{i=1}^2 |\nabla \mu_i\om - \nabla \mu_0\om|^2\right)(x,t)\dd x\dd t \le C^\ve.
\end{align}
Furthermore, differentiating \eqref{e3m} with respect to the spatial variables we obtain that
\begin{align}
\nonumber
&\nu {\nabla\dot w}\om + E(\vp_1\om,\vp_2\om)\, \nabla w\om + w\om \left(\partial_1 E(\vp_1\om,\vp_2\om)\nabla\vp_1\om +\partial_2 E(\vp_1\om,\vp_2\om)\nabla\vp_2\om\right) \\ \label{es11a}
&\quad + (\nabla w\om {-} \nabla\vp_0\om) \frac{(f^{-1})'(\vp_0\om {-} w\om)}{|\vp_0\om|{+}|\vp_1\om|{+}|\vp_2\om|}
=  f^{-1}(\vp_0\om {-} w\om)\nabla\left(\frac{1}{|\vp_0\om|{+}|\vp_1\om|{+}|\vp_2\om|}\right).
\end{align}
Testing \eqref{es11a} by $\nabla w\om$ and using \red{the inequality $|\vp_0\om|{+}|\vp_1\om|{+}|\vp_2\om| \ge 1$, Hypotheses \ref{h1}\,(ii), (iv), and the estimate  \eqref{es10}, we find that}
\be{es11}
\frac{\dd}{\dd t}|\nabla{w}\om|_H^2+ c\, |\nabla{w}\om|_H^2
\le C\left(1+ |\nabla\vp_0\om|_H^2 + |\nabla\vp_1\om|_H^2 + |\nabla\vp_2\om|_H^2\right)
\ee
with some constants \vio{$C>c \ge 0$.} From \eqref{es0} we immediately obtain the pointwise bound
\be{es12}
|\nabla\vp_0\om| \le |\nabla\vp_1\om|+|\nabla\vp_2\om|\quad \mbox{a.\,e.}
\ee
It follows from \eqref{es9a}, \eqref{es11}, and by comparison in \eqref{es11a} that
\be{es11b}
\supess_{t\in (0,T_m)} \left(|{\nabla\dot w}\om|_H + \vio{|{\nabla w}\om|_H}\right) \le C.
\ee
By virtue of \eqref{e4m} we have
\be{es13}
|\nabla\mu_0\om| \le |\nabla\vp_1\om|+|\nabla\vp_2\om|+|\nabla w\om|\quad \mbox{a.\,e.}
\ee
Since $\nabla_{\!\vp}\psi^\ve$ is Lipschitz continuous for every $\ve>0$, we obtain from \cyan{\eqref{e4mm} and \eqref{es10}} that
\be{es14}
\red{|\bar\mu_i\om(t)|} \le C^\ve \left(1+ \io\sum_{i=1}^2|\vp_i\om|^2(x,t)\dd x\right)^{1/2}
\ee
and
\be{es15}
\io|\mu_0\om|^2\dd x \le C\left(1 + \io|\vp_0\om|^2\dd x\right). 
\ee
We now summarize the above computations in \eqref{es9a}--\eqref{es15} and obtain for all $t \in (0,T_m)$ that
\begin{align}\nonumber
&\io\left(\sum_{i=0}^2 \left(|\vp_i\om|^2  + |\nabla \vp_i\om|^2\right) + |\vr\om|^2  + |\nabla \vr\om|^2 + |\nabla w\om|^2 + |\nabla \dot w\om|^2\right)(x,t)\dd x\\ \label{es16}
& \hspace{10mm} + \int_0^t\io \left(\sum_{i=0}^2 \left(|\mu_i\om|^2 + |\nabla \mu_i\om|^2\right) +|\dot \vr\om|^2\right)(x,\tau) \dd x\dd \tau \le C^\ve
\end{align}
with a constant $C^\ve>0$ independent of $m$, \cyan{and the uniform estimate \eqref{es10} holds.} By comparison in \eqref{e4mm}, we have a bound for $\Delta \vp_i\om$ in $L^2(\Omega\times (0,T))$ independent of $m$, $i=1,2$. Finally, by comparison in \eqref{e1m}, we obtain bounds in $L^2(0,T;W^{-1,2}(\Omega))$ independent of $m$ for $\dot\vp_i\om$, $i=0,1,2$. We thus have sufficient estimates which on the one hand guarantee that the solution exists on the whole time interval $[0,T]$ and, on the other hand, enable us to pass to the limit as $m \to \infty$ in \eqref{e1m}--\eqref{qrm} and check that the following statement holds.
\begin{propo}\label{p1}
Let Hypothesis \ref{h1} hold and let $\ve>0$ be given. Then the system \eqref{e1w}--\eqref{e7w}, \eqref{ini1} admits a solution with the regularity $\mu_i, \nabla \mu_i, \Delta \vp_i \in L^2((0,T)\times\Omega)$, $\red{\vp_i,} \nabla \vp_i \in L^\infty(0,T;L^2(\Omega))$, $\dot \vp_i \in L^2(0,T;W^{-1,2}(\Omega))$ for $i=0,1,2$, $\vp_0+\vp_1+ \vp_2 =1$ a.\,e., $w,\cyan{\dot w} \in L^\infty((0,T)\times\Omega)$, $\nabla w, \nabla\dot w \in L^\infty(0,T;L^2(\Omega))$, $\dot\vr \in L^2((0,T)\times\Omega)$, $\vr, \nabla \vr \in L^\infty(0,T;L^2(\Omega))$.
\end{propo}

We can indeed pass to the limit in the initial conditions for $\vr$ and $w$ by virtue of \eqref{esro3} and \eqref{es11b}. For the initial conditions for $\vp_i$, the argument is standard as well. It is easy to check for each $i=0,1,2$ that
\be{weakini}
\begin{aligned}
&\forall \eta > 0\ \forall v_i \in L^2(\Omega) \ \exists t_\eta > 0:\\
&\qquad t\in (0,t_\eta)\ \Longrightarrow \ \exists m_\eta \in \nat \ \forall m>m_\eta: \ \left|\io(\vp_i\om(x,t) - \vp_i\om(x,0))v_i(x)\dd x\right| < \eta,
\end{aligned}
\ee
so that the initial condition is satisfied in weak sense. \blu{Note that this is related to the so-called Aubin-Lions lemma, see \cite[Theorem 5.1]{lions} for the original reference.}


\section{Limit as $\ve \to 0$}\label{lim}

In the previous section we have proved that system \eqref{e1w}--\eqref{e7w}, \eqref{ini1} admits a global solution. The estimates that we have derived so far depend on $\ve$. We split this section into two subsections. In Subsection \ref{inde}, we derive estimates independent of $\ve$ of the solution of \eqref{e1w}--\eqref{e7w}, and in Subsection \ref{prt1} we prove Theorem \ref{t1} by passing to the limit as $\ve \to 0$.


\subsection{Estimates independent of $\ve$}\label{inde}

Let us start with the following simple modification of \cite[Propositions 2.10, 2.13]{ckrs}.

\begin{propo}\label{ckrs}
Let $\psi$ satisfy Hypothesis \ref{h1}\,(v). Then there exists
$\bar\ve>0$ and positive constants $b,c,r$ such that for $\ve \in (0,\bar\ve)$, the Yosida approximations $\psi^\ve$ of $\psi$ have the following properties:
\begin{itemize}
\item[{\rm (i)}] $\dist(\hat\vp,\Theta_{\delta_T}) \le \delta_T/2 \ \Longrightarrow \ |\nabla_{\!\vp}\psi^\ve(\hat\vp)| \le b$;
\item[{\rm (ii)}] $\hat\vp \in \Theta_{\delta_T},\ \vp\in \real^2,\ |\vp-\hat\vp|\ge \delta_T/2$
\item[] \qquad $\Longrightarrow\ r|\nabla_{\!\vp}\psi^\ve(\vp)-\nabla_{\!\vp}\psi^\ve(\hat\vp)|\leq \scal{\nabla_{\!\vp}\psi^\ve(\vp)-\nabla_{\!\vp}\psi^\ve(\hat\vp),\vp - \hat\vp}+ c$.
\end{itemize}
\end{propo}

\bpf{Proof}
We prove the statement for $b=b'$, $r = r'$, $c=c'+2r'b'$, where $b',c',r'$ are as in Hypothesis \ref{h1}\,(v).
Let us start with part (i), and consider $\hat \vp \in \real^2$ such that $\dist(\hat\vp,\Theta_{\delta_T}) \le \delta_T/2$. For $\ve>0$ we define $J^\ve \hat\vp$ as in Proposition \ref{yos}, and choose any $\hat\xi \in \partial\psi(\hat\vp)$. We have by \eqref{yos1} that
$$
\hat\xi^\ve:= \nabla_{\!\vp}\psi^\ve(\hat\vp) = \frac{1}{\ve}(\hat\vp-J^\ve \hat\vp) \in \partial\psi(J^\ve \hat\vp),
$$
hence \red{$-\ve\scal{\hat\xi^\ve - \hat\xi, \hat\xi^\ve} = \scal{\hat\xi^\ve - \hat\xi, J^\ve \hat\vp-\hat\vp} \ge 0$} by the monotonicity of $\partial\psi$. We thus have $|\hat\xi^\ve| \le b'$ by Hypothesis \ref{h1}\,(v), and part (i) is proved.
\\[2mm]
To prove part (ii), let $\hat \vp\in\Theta_{\delta_T}$ be given, and put $\bar\ve = \delta_T/(4b')$. For $\ve < \bar \ve$ we have
$$
|\hat\vp - J^\ve \hat\vp| = \ve |\hat\xi^\ve| < \frac{\delta_T}{4} 
$$
by virtue of part (i). Hence, $\dist(J^\ve \hat\vp,\Theta_{\delta_T}) < \delta_T/4$. Let further $|\vp-\hat\vp|\ge \delta_T/2$ for some $\vp \in \real^2$. We denote $\xi^\ve = \nabla_{\!\vp}\psi^\ve(\vp)$.
We have either 
\be{y1}
|J^\ve \vp - J^\ve \hat\vp| < \frac{\delta_T}{4},
\ee
or
\be{y2}
|J^\ve \vp - J^\ve \hat\vp| \ge \frac{\delta_T}{4}.
\ee
In case \eqref{y1}, we have $\dist(J^\ve \hat\vp,\Theta_{\delta_T}) < \delta_T/2$, and we obtain from (i) simply that
$$
r|\xi^\ve-\hat\xi^\ve|\le r(|\xi^\ve| +|\hat\xi^\ve|)  \le 2rb.
$$
If \eqref{y2} holds, then we have by Hypothesis \ref{h1}\,(v) that
$$
r|\xi^\ve-\hat\xi^\ve|\le \scal{\xi^\ve{-}\hat\xi^\ve,J^\ve \vp {-} J^\ve \hat\vp} + c' = \scal{\xi^\ve{-}\hat\xi^\ve,\vp {-} \hat\vp} - \ve|\xi^\ve{-}\hat\xi^\ve|^2 +c' \le\scal{\xi^\ve{-}\hat\xi^\ve,\vp {-} \hat\vp} +c'.
$$
Combining the two inequalities and using the monotonicity of $\nabla_{\!\vp}\psi^\ve$ we obtain the assertion.
\epf

We actually need the following consequence of Proposition \ref{ckrs}.

\begin{coro}\label{c1}
Let $\psi, \bar\ve, b, c, r$ be as in Proposition \ref{ckrs}. Then there exists a constant $\hat c>0$ with the property that for every $\ve < \bar\ve$, for every $\hat\vp \in L^2(\Omega)$ such that 
$\hat \vp(x) \in \Theta_{\delta_T}$ a.\,e., and for every $\vp \in L^2(\Omega)$ we have
\be{cor1}
r \io |\nabla_{\!\vp}\psi^\ve(\vp(x))-\nabla_{\!\vp}\psi^\ve(\hat\vp(x))|\dd x \le
\io \scal{\nabla_{\!\vp}\psi^\ve(\vp(x))-\nabla_{\!\vp}\psi^\ve(\hat\vp(x)),\vp(x) - \hat\vp(x)}\dd x + \hat c.
\ee
\end{coro}

\bpf{Proof}
Let $\vp \in L^2(\Omega)$ be arbitrarily chosen. We define $\Omega_+ := \{x \in \Omega: \dist(\vp(x), \Theta_{\delta_T}) \ge \delta_T/4\}$, $\Omega_- = \Omega \setminus \Omega_+$. For a.\,e. $x \in \Omega_-$ we have by Proposition \ref{ckrs}\,(i) that
$$
|\nabla_{\!\vp}\psi^\ve(\vp(x))-\nabla_{\!\vp}\psi^\ve(\hat\vp(x))| \le 2b\,.
$$
For a.\,e. $x \in \Omega_+$ Proposition \ref{ckrs}\,(ii) yields that
$$
r |\nabla_{\!\vp}\psi^\ve(\vp(x))-\nabla_{\!\vp}\psi^\ve(\hat\vp(x))| \le
\scal{\nabla_{\!\vp}\psi^\ve(\vp(x))-\nabla_{\!\vp}\psi^\ve(\hat\vp(x)),\vp(x) - \hat\vp(x)}+ c\,.
$$
Using the fact that $\scal{\nabla_{\!\vp}\psi^\ve(\vp(x))-\nabla_{\!\vp}\psi^\ve(\hat\vp(x)),\vp(x) - \hat\vp(x)} \ge 0$ a.\,e. we can combine the two inequalities and obtain that
\begin{align*}
r \io |\nabla_{\!\vp}\psi^\ve(\vp(x))-\nabla_{\!\vp}\psi^\ve(\hat\vp(x))|\dd x &\le
\io \scal{\nabla_{\!\vp}\psi^\ve(\vp(x))-\nabla_{\!\vp}\psi^\ve(\hat\vp(x)),\vp(x) - \hat\vp(x)}\dd r \\
& \ + c|\Omega_+| + 2rb|\Omega_-|\,.
\end{align*}
Putting $\hat c := |\Omega|(c+2rb)$ we complete the proof.
\epf

\blu{Let us come back to Problem \eqref{e1w}--\eqref{e7w} with initial conditions \eqref{ini1} and} estimate the distance of the functions $\bar\vp_i(t)$ from the boundary of $\Theta$. To this end, we choose $v_0 = v_1 = v_2 = 1$ and put
$$
\Gamma =\frac{\gamma(\vr)}{(|\vp_0|+|\vp_1|+|\vp_2|)(|\bar\vp_0|+|\bar\vp_1|+|\bar\vp_2|)},
$$
We obtain
\begin{align}\label{pe1}
\dot{\bar\vp}_0(t) &= \blu{-\frac{\bar\vp_0(t)}{|\Omega|}}\io \Gamma(x,t)\,(1-\vp_0(x,t))\dd x,\\ \label{pe2}
\dot{\bar\vp}_1(t) &= \frac{\bar\vp_0(t)}{|\Omega|}\io \Gamma(x,t)\,\vp_1(x,t)\dd x,\\ \label{pe3}
\dot{\bar\vp}_2(t) &= \frac{\bar\vp_0(t)}{|\Omega|}\io \Gamma(x,t)\,\vp_2(x,t)\dd x.
\end{align}
From Hypothesis \ref{h1}\,(iii) it follows that $|\Gamma(x,t)\,(1-\vp_0(x,t))| \le |\Gamma(x,t)|\,(|\vp_1(x,t)|+|\vp_2(x,t)|) \le K$ for a.\,e. $(x,t) \in \Omega\times (0,T)$. By Hypothesis \ref{h1}\,(vii) we have $\bar\vp_0(0) \ge \delta/\sqrt{2} > 0$, hence,
\be{pe4}
\bar\vp_0(t) \ge  \bar\vp_0(0)\expe^{-Kt} > 0 \ \ \mbox{for all } \ t \in [0,T]. 
\ee
Lower bounds for $\bar\vp_1, \bar\vp_2$ are more delicate to obtain.
These functions are continuously differentiable. There exists therefore $T_\ve \in [0,T]$ such that
\be{pe5}
\bar\vp_i(t) \ge \delta \expe^{-KT-1} \ \ \mbox{for all } \ t \in [0,T_\ve], \ i=1,2. 
\ee
Put $T_\ve^* = \max\{T_\ve \in [0,T]: \mbox{ inequality \eqref{pe5} holds}\}$, and assume that $T_\ve^*<T$ for some $\ve<\bar\ve$. For definiteness, we can assume that
\be{pe6}
\bar\vp_1(T_\ve^*) = \delta \expe^{-KT-1}. 
\ee
Taking into account \eqref{pe4}, we have $1-\bar \vp_1(t) - \bar\vp_2(t) = \bar\vp_0(t) > (\delta/2)\expe^{-Kt}$ in $[0, T_\ve^*]$. Hence, denoting $\vp=(\vp_1, \vp_2)$ we have $\dist(\bar\vp(t), \partial\Theta) \ge (\delta/2)\expe^{-KT-1} > \delta_T$, so that
$\bar\vp(t) \in \Theta_{\delta_T}$ for all $t \in [0,T_\ve^*]$.

\cyan{Recall that we have the bound
\be{es10a}
\supess_{(x,t)\in \Omega\times (0,T_\ve^*)} (|w(x,t)| \cyan{+ |\dot w(x,t)|}) \le C
\ee
as a consequence of \eqref{es10}.} Let us denote $\bar\mu = (\bar\mu_1, \bar\mu_2)$.  From \cyan{\eqref{es10a}}, \eqref{e4w}, and \eqref{cg} it follows that
\be{pe7}
|\bar\mu(t)| \le \io\left(|\nabla_{\!\vp}\psi^\ve(\vp)|+ |\nabla_{\!\vp}g(\vp)|+ \cyan{\frac12|\nabla_{\!\vp}E(\vp)| w^2}\right))\dd x \le \io|\nabla_{\!\vp}\psi^\ve(\vp){-}\nabla_{\!\vp}\psi^\ve(\bar\vp)|\dd x + \cyan{\tilde C}|\Omega|
\ee
with \cyan{$\tilde C= b+C_g + \frac12\sup|\nabla_{\!\vp}E(\vp)|  w^2$} 
 where we have used Hypothesis \ref{h1}\,(ii),(vi), and Proposition \ref{ckrs}\,(i). We further obtain from Corollary \ref{c1} and \eqref{e4w} that
\begin{align}\nonumber
|\bar\mu(t)| &\le \frac{1}{r}\io \scal{\nabla_{\!\vp}\psi^\ve(\vp)-\nabla_{\!\vp}\psi^\ve(\bar\vp),\vp - \bar\vp}\dd x + \cyan{\tilde C}|\Omega| + \frac{c}{r}\\ \nonumber
&= \frac{1}{r}\io\scal{\nabla_{\!\vp}\psi^\ve(\vp),\vp - \bar\vp}\dd x + \cyan{\tilde C}|\Omega| + \frac{c}{r}\\ \label{pe8}
&= \frac{1}{r}\left(-\io |\nabla\vp|^2\dd x - \io\scal{\nabla_{\!\vp}g(\vp), \vp - \bar \vp}\dd x + \io\scal{\mu, \vp - \bar\vp}\dd x \right) + \cyan{\tilde C}|\Omega| + \frac{c}{r}\,.
\end{align}
We now use again \eqref{cg}, the fact that 
$\io\scal{\mu, \vp - \bar\vp}\dd x = \io\scal{\mu-\bar \mu, \vp - \bar\vp}\dd x$, and the elementary inequalities
\be{ele}
\io|\vp - \bar\vp|^2\dd x \le C \io|\nabla\vp|^2\dd x, \ \ 
\io|\mu - \bar\mu|^2\dd x \le C \io|\nabla\mu|^2\dd x.
\ee
to conclude that there exists a constant M independent of $\ve$ such that for all $t \in [0,T^*_\ve]$ we have
\be{pe9}
|\bar\mu(t)| \le M\left(1+\left(\io |\nabla\vp|^2(x,t)\dd x\right)^{1/2}\left(\io |\nabla\mu|^2(x,t)\dd x\right)^{1/2}\right).
\ee
We now repeat the estimation procedure from Subsection \ref{es01}. We test the $i$-th equation in \eqref{e1w} by \vio{$v_i=\mu_i$,} and sum up to obtain similarly as in \eqref{es1}--\eqref{es4}
\red{
\begin{align}\nonumber
&\frac{\dd}{\dd t}\io\mathcal{F}_0^\ve(\vp_0, \vp_1, \vp_2, w)\dd x + \frac{c}{2}\io \sum_{i=1}^{2}(|\nabla\mu_i - \nabla \mu_0|^2)\dd x\\ \label{pe10}
&\hspace{10mm}\le \sum_{i=0}^2\io S_i\mu_i\dd x + \cyan{C\left(1+ \io|\vp_0|\dd x\right)}
\end{align} 
for a.\,e. $t \in (0,T_\ve^*)$} \cyan{with $\mathcal{F}_0^\ve$ defined in \eqref{poteps}, and with some constants $C>c>0$ independent of $\ve$. We further estimate the right-hand side by}
$$
\sum_{i=0}^2\io S_i\mu_i\dd x \le C\sum_{i=0}^2\io |\vp_i||\mu_i|\dd x,
$$
to obtain, by virtue of \eqref{ele}, \eqref{pe9}, and the hypotheses on $\hat F$ and $g$, that
\begin{align}\nonumber
&\io\left(|\vp_0|^2 + \psi^\ve(\vp) +|\nabla\vp|^2\right)(x,t)\dd x + \int_0^t\io \sum_{i=1}^2(|\nabla\mu_i - \nabla\mu_0|^2)(x,\tau)\dd x\dd\tau\\ \label{pe11}
&\hspace{3mm}\le C\left(1 + \int_0^t\io\left(|\nabla\vp|^2+|\vp_0|^2 + |\nabla\mu_1|^2 + |\nabla\mu_2|^2\right)(x,t)\dd x\dd\tau\right).
\end{align} 
We have for all $(x,t) \in \Omega\times [0,T^*_\ve]$ the identity $\vp_0 + \vp_1 + \vp_2 = 1$ and $\nabla\vp_0 = - \nabla\vp_1 - \nabla\vp_2$, hence, $|\nabla\mu_0| \le C(|\nabla\vp_0| + |\nabla w|)$. Note that
repeating the computations leading to \eqref{es11b} we derive the estimates
\be{pe12a}
 \supess_{t\in (0,T_\ve^*)} \left(|{\nabla\dot w}|_H +|{\nabla w}|_H\right) \le C
\ee
with $C$ independent of $\ve$. The Gronwall argument and \eqref{lowps} now yield
\be{pe12}
\io\left(\psi^\ve(\vp) + \sum_{i=0}^2\left(|\vp_i|^2 + |\nabla\vp_i|^2\right)\right)(x,t)\dd x + \int_0^t\io \sum_{i=1}^2\left(|\mu_i|^2 +|\nabla\mu_i|^2\right) (x,\tau)\dd x\dd\tau \le C^*
\ee
for every $t \in [0,T_\ve^*]$ with a constant $C^*>0$ independent of $\ve$. By comparison in \eqref{e1w} we get the bound
\be{pe12b}
\int_0^{T_\ve^*} \|\dot \vp_i(t)\|_{W^{-1,2}(\Omega)}^2 \dd t \le C\,, \quad i=0,1,2\,.
\ee
To make the list of estimates complete, recall that the upper bound in \eqref{esro}--\eqref{esro2} is independent of $m$ and $\ve$, so that
\be{pe12c}
\supess_{t\in (0,T_\ve^*)} \left(|\vr(t)|_H +|\nabla\vr(t)|_H\right) \le C, \quad \int_0^{T_\ve^*} |\dot\vr(t)|_H^2 \dd t \le C\,.
\ee
The next step consists in proving that $T^*_\ve = T$. To this end, we split for each $t \in [0,T_\ve^*]$ the domain $\Omega$ into three parts, namely
\begin{align*}
\Omega_0(t) &= \{x \in \Omega: \vp_1(x,t) \ge 0\}, \\
\Omega_1(t) &= \{x \in \Omega: 0 > \vp_1(x,t) \ge -\ve^{1/4}\},\\
\Omega_2(t) &= \{x \in \Omega: -\ve^{1/4} > \vp_1(x,t)\}.
\end{align*}
Let us start with $\Omega_2(t)$. By definition \eqref{yosi} of $\psi^\ve$, we have for $x \in \Omega_2(t)$ that
\be{ome2a}
\psi^\ve(\vp(x,t)) \ge \frac{1}{2\ve}\min_{z \in \Theta}|\vp_1(x,t) - z_1|^2 \ge \frac{1}{2\sqrt{\ve}}.
\ee
By virtue of \eqref{pe12}, we have
\be{ome2b}
|\Omega_2(t)| \le \red{2 C^*\sqrt{\ve}.}
\ee
We now rewrite Eq.~\eqref{pe2} in the form
$$
\dot{\bar\vp}_1(t) = \vio{\frac{\bar\vp_0(t)}{|\Omega|}}\left(\int_{\Omega_0(t)} + \int_{\Omega_1(t)}+ \int_{\Omega_2(t)}\right)\Gamma(x,t)\,\vp_1(x,t)\dd x, 
$$
where \red{
\begin{align*}
\int_{\Omega_1(t)}\Gamma(x,t)\,\vp_1(x,t)\dd x &\ge -K|\Omega|\ve^{1/4},\\
\int_{\Omega_2(t)}\Gamma(x,t)\,\vp_1(x,t)\dd x &\ge -K\int_{\Omega_2(t)}|\vp_1(x,t)|\dd x\ge -K|\Omega_2(t)|^{1/2}\left(\io|\vp_1(x,t)|^2\dd x\right)^{1/2}\\
&  \ge -\sqrt{2}KC^*\ve^{1/4},\\
\int_{\Omega_0(t)}\Gamma(x,t)\,\vp_1(x,t)\dd x &\ge -K\int_{\Omega_0(t)}\vp_1(x,t)\dd x\\
& = -K |\Omega| \bar\vp_1(t) + K\left(\int_{\Omega_1(t)}{+} \int_{\Omega_2(t)}\right) \vp_1(x,t)\dd x\\
&\ge -K |\Omega| \bar\vp_1(t) - K(1+\sqrt{2}C^*)\, \ve^{1/4}.
\end{align*}
Using the fact that $0\le \bar\vp_0(t) \le 1$ for $t \in [0,T_\ve^*]$ we have that}
\be{pe13}
\dot{\bar\vp}_1(t) = \red{\frac{\bar\vp_0(t)}{|\Omega|}}\int_{\Omega}\Gamma(x,t)\,\vp_1(x,t)\dd x \ge -K(\bar\vp_1(t) +\Lambda\ve^{1/4})
\ee
with a constant $\Lambda>0$ independent of $\ve$. We thus obtain a lower bound for $\bar\vp_1(t)$, namely (note that $\bar\vp_1(0) \ge \delta$ by Hypothesis \ref{h1}\,(vii)),
\be{pe14}
\bar\vp_1(t) \ge  \delta \expe^{-Kt} -\Lambda\ve^{1/4}\left(1 - \expe^{-Kt}\right) \ge \delta \expe^{-Kt} -\Lambda\ve^{1/4}
\ee
for $t \in [0,T_\ve^*]$. We see that for $\ve>0$ sufficiently small, condition \eqref{pe6} is violated. Hence, by \eqref{pe5}, $T^*_\ve =T$ and the estimate \eqref{pe12} holds globally in $[0,T]$.


\subsection{Proof of Theorem \ref{t1}}\label{prt1}

We show that passing to the limit as $\ve \to 0$ in \eqref{e1w}--\eqref{e7w} we obtain a solution to \eqref{e1}--\eqref{e7} in the sense of Theorem \ref{t1}. We label here the solution $(\mu_i, \vp_i, w, \vr)$ of \eqref{e1w}--\eqref{e7w} with the upper index $\ve$ in order to emphasize the dependence on $\ve$.

The estimates \eqref{pe12}--\eqref{pe12c} are independent of $\ve$ and hold globally on $[0,T]$. We can therefore extract a subsequence $\ve \to 0$ such that
\begin{itemize}
\item $\nabla\vp_i^\ve \to \nabla\vp_i$ for $i=0,1,2$, $\nabla \vr^\ve \to \nabla \vr$, $\nabla w^\ve \to \nabla w$ weakly-star in $L^\infty(0,T; L^2(\Omega))$;
\item $\dot\vr^\ve \to \dot\vr$, $\mu_i^\ve \to \mu$, $\nabla\mu_i^\ve \to \nabla \mu_i$ for $i=0,1,2$, $\dot w^\ve \to \dot w$ weakly in $L^2(\Omega\times (0,T))$;
\item $\dot\vp_i^\ve \to \dot\vp_i$ for $i=0,1,2$ weakly in $L^2(0,T; W^{-1,2}(\Omega))$.
\end{itemize}
Using the Sobolev embedding theorems, the trace theorem, and the Lions compactness lemma \cite[Theorem 5.1]{lions} we obtain the convergences, passing again to a subsequence of $\ve \to 0$ if necessary,
\begin{itemize}
\item $\vr^\ve \to \vr$, $w^\ve \to w$ strongly in $C([0,T]; L^2(\Omega))$;
\item $\vp_i^\ve \to \vp_i$ for $i=0,1,2$ strongly in $L^2(\Omega\times (0,T))$;
\item $\vr^\ve \to \vr$ strongly in $L^2(0,T; L^2(\partial\Omega))$.
\end{itemize}
We can pass to the limit in all terms in \eqref{e1w}--\eqref{e7w}, and the limit initial condition \eqref{ini1} is obtained by an argument similar to \eqref{weakini}. The variational inequality \eqref{e4u} needs to be paid some attention. Since $\psi^\ve$ is convex, we can rewrite \eqref{e4w} as 
\begin{align}\nonumber
&\io\left(\mu^\ve_1 - \partial_1 g(\vp^\ve_1, \vp^\ve_2) \cyan{-\partial_1 E(\vp^\ve_1,\vp^\ve_2) \frac{w^2}{2}}\right)(v_1 - \vp^\ve_1)\dd x\\ \nonumber
&\quad + \io\left(\mu^\ve_2 - \partial_2 g(\vp^\ve_1, \vp^\ve_2)\cyan{-\partial_2 E(\vp^\ve_1,\vp^\ve_2) \frac{w^2}{2}}\right)(v_2 - \vp^\ve_2)\dd x\\ \label{e4e}
&\quad -\io\big(\scal{\nabla\vp^\ve_1,\nabla(v_1{-}\vp^\ve_1)} + \scal{\nabla\vp^\ve_2,\nabla(v_2{-}\vp^\ve_2)}\big)\dd x \le \io\big(\psi^\ve(v_1, v_2) -\psi^\ve(\vp^\ve_1, \vp^\ve_2) \big)\dd x
\end{align}
for a.\,e. $t \in (0,T)$ and for all test functions $v_1, v_2\in W^{1,2}(\Omega)$. We now choose an arbitrary test function $\lambda \in L^2(0,T)$, $\lambda(t) \ge 0$ a.\,e. From the above convergences it follows that 
$$
\liminf_{\ve \to 0} \int_0^T\io|\nabla \vp^\ve_i(x,t)|^2\lambda(t)\dd x\dd t \ge \int_0^T\io|\nabla \vp_i(x,t)|^2\lambda(t)\dd x\dd t\,,
$$
and using \eqref{yos4} we obtain the pointwise limit $\lim_{\ve \to 0} \psi^\ve(v_1, v_2) =\psi(v_1, v_2)$. We multiply both sides of the inequality \eqref{e4e} by $\lambda(t)$, integrate over $t\in (0,T)$ and pass to the limit to obtain
\begin{align}\nonumber
&\int_0^T\io\left(\mu_1 - \partial_1 g(\vp_1, \vp_2)\cyan{-\partial_1 E(\vp_1,\vp_2) \frac{w^2}{2}}\right)(v_1 -\vp_1)\lambda(t)\dd x\dd t \nonumber\\
&\qquad + \int_0^T\io\left(\mu_2 - \partial_2 g(\vp_1, \vp_2)\cyan{-\partial_2 E(\vp_1,\vp_2) \frac{w^2}{2}}\right)(v_2-\vp_2)\lambda(t)\dd x\dd t \nonumber\\
&\qquad -\int_0^T\io\big(\scal{\nabla\vp_1,\nabla(v_1{-}\vp_1)} + \scal{\nabla\vp_2,\nabla(v_2{-}\vp_2)}\big)\lambda(t)\dd x\dd t\nonumber\\ \label{e4l}
&\quad\le \int_0^T\io\psi(v_1, v_2)\lambda(t)\dd x\dd t -\liminf_{\ve \to 0} \int_0^T\io\psi^\ve(\vp^\ve_1, \vp^\ve_2)\lambda(t)\dd x\dd t
\end{align}
for all test functions $v_1, v_2\in W^{1,2}(\Omega)$. It remains to prove that we have
\be{psil}
\liminf_{\ve \to 0} \int_0^T\io \psi^\ve(\vp^\ve_1(x,t), \vp^\ve_2(x,t))\lambda(t)\dd x\dd t\ge \int_0^T\io \psi(\vp_1(x,t), \vp_2(x,t))\lambda(t)\dd x\dd t.
\ee
If \eqref{psil} is fulfilled, then, on the one hand, \eqref{e4u} holds and, on the other hand, we conclude that $\psi(\vp_1(x,t), \vp_2(x,t)) < \infty$ almost everywhere. This means in particular that $(\vp_1(x,t), \vp_2(x,t)) \in \Theta$ for a.\,e. $(x,t) \in \Omega\times (0,T)$. Hence, as mentioned on the last line of Section \ref{sta}, the identity $|\vp_0|+|\vp_1|+|\vp_2| = \vp_0+\vp_1+\vp_2 = 1$ holds almost everywhere, so that \eqref{e3w} coincides with \eqref{e3}, and \eqref{e5w}--\eqref{qr} coincides with \eqref{e5}.

To prove \eqref{psil}, we first notice that by \eqref{pe10} we have
$$
\supess_{t \in (0,T)} \io \psi^\ve(\vp^\ve(x,t))\dd x \le C.
$$
For simplicity, we omit for a moment the arguments $(x,t)$ and write simply $\vp^\ve, \vp$ instead of $\vp^\ve(x,t), \vp(x,t)$.
By \eqref{yos5}, we have
\be{el4}
\psi^\ve(\vp^\ve) \ge \frac{1}{2\ve} |\vp^\ve - J^\ve\vp^\ve|^2 \quad \mbox{a.\,e.}
\ee
Hence, for a.\,e. $t\ge 0$,
\be{el5}
\io |\vp^\ve - J^\ve\vp^\ve|^2 \dd x \le 2\ve \io \psi^\ve(\vp^\ve)\dd x \le C\ve\,.
\ee
We thus have {for a.\,e. $t\ge 0$} by triangle inequality
\be{el1}
|J^\ve\vp^\ve(t) - \vp(t)|_H \le |J^\ve\vp^\ve(t) - \vp^\ve(t)|_H + |\vp^\ve(t) - \vp(t)|_H \le C\ve + |\vp^\ve(t) - \vp(t)|_H.
\ee
We know that $\vp^\ve$ converge to $\vp$ in $L^2(\Omega\times(0,T))$. In particular, it follows from \eqref{el1} that $J^\ve\vp^\ve(x,t) \to \vp(x,t)$ a.\,e. in $\Omega\times(0,T)$. On the other hand, by \eqref{yos3} we have
\be{el2}
\psi^\ve(\vp^\ve) \ge \psi(J^\ve\vp^\ve) \quad \mbox{a.\,e.},
\ee
and \eqref{psil} follows from \eqref{el1}--\eqref{el2} and from the lower semicontinuity of $\psi$. We thus obtain the inequality
\begin{align}\nonumber
&\int_0^T\io\left(\left(\mu_1 - \partial_1 g(\vp_1, \vp_2)\cyan{-\partial_1 E(\vp_1,\vp_2) \frac{w^2}{2}}\right)(v_1{-}\vp_1)\right. \\
\nonumber 
&\qquad \left.+ \left(\mu_2 - \partial_2 g(\vp_1, \vp_2)\cyan{-\partial_2 E(\vp_1,\vp_2) \frac{w^2}{2}}\right)(v_2{-}\vp_2)\right)\lambda(t)\dd x\dd t\\ \nonumber
&\qquad -\int_0^T\io\big(\scal{\nabla\vp_1,\nabla(v_1{-}\vp_1)} + \scal{\nabla\vp_2,\nabla(v_2{-}\vp_2)}\big)\lambda(t)\dd x\dd t\\ \label{e4z}
&\quad\le \int_0^T\io\psi(v_1, v_2)\lambda(t)\dd x\dd t -\int_0^T\io\psi(\vp_1, \vp_2)\lambda(t)\dd x\dd t
\end{align}
for all test functions $v_1, v_2\in W^{1,2}(\Omega)$, $\lambda \in L^2(0,T)$, $\lambda(t) \ge 0$ a.\,e., which is equivalent to \eqref{e4u}.
This completes the proof of Theorem \ref{t1}.


\section*{Acknowledgments}

This research was supported by the Italian Ministry of Education, University and Research (MIUR): 
Dipartimenti di Eccellenza Program (2018--2022) -- Dept. of Mathematics ``F. Casorati'', University of Pavia, by the GA\v CR Grant No.~20-14736S and by the European Regional Development Fund, Project No. CZ.02.1.01/0.0/0.0/16{\_}019/0000778.
In addition, it has been performed in the framework of the project Fondazione Cariplo-Regione Lombardia  MEGAsTAR ``Matema\-tica d'Eccellenza in biologia ed ingegneria come acceleratore di una nuova strateGia per l'ATtRattivit\`a dell'ateneo pavese''.
The present paper also benefits from the support of the  GNAMPA (Gruppo Nazionale per l'Analisi Matematica, la Probabilit\`a e le loro Applicazioni)
of INdAM (Istituto Nazionale di Alta Matematica) for E.\,R.


\end{document}